\documentclass[11pt]{amsart}
\usepackage{amsmath,amsthm,amsfonts,amssymb}

\newcommand{\C}{{\mathbf C}}
\newcommand{\N}{{\mathbf N}}

\newcommand{\Q}{{\mathbf Q}}

\renewcommand{\P}{{\mathbf P}}
\newcommand{\R}{{\mathbf R}}

\newcommand{\Sing}{\mathrm{ Sing}}

\newcommand{\ortho}[1]{{#1}^{\bot}}
\newcommand{\st}{\mathrm{Stab}}

\newtheorem{remark}[equation]{Remark}

\newtheorem{remarks}[equation]{Remarks}
\newtheorem{lemma}[equation]{Lemma}
\newtheorem{proposition}[equation]{Proposition}
\newtheorem{theorem}[equation]{Theorem}
\numberwithin{equation}{section}

\title{%
A note on Jets of Entire Curves in Semi-Abelian Varieties
}
\author {Junjiro Noguchi and J\"org Winkelmann
}
\address{
Noguchi, Junjiro\\
Graduate School of Mathematical Sciences\\
University of Tokyo\\
Komaba, Meguro,Tokyo 153-8914\\
Japan
}
\email{noguchi@ms.u-tokyo.ac.jp
}
\address{%
Winkelmann, J\"org \\
Korean Institute for Advanced Studies\\
School of Mathematics\\
207-43 Cheongryangri-dong\\
Dongdaemun-gu\\
Seoul\\
130-012 Korea
}
\email{jwinkel@member.ams.org\newline\indent{\itshape Webpage: }%
http://www.math.unibas.ch/\~{ }winkel/
}
\begin{document} 
\baselineskip18pt
\maketitle
\thispagestyle{empty}
\footnotetext{  Research supported in part by Grant-in-Aid
for Scientific Research (A)(1), 13304009.}


\begin{abstract}
We prove a product decomposition of the Zariski closure of
the jet lifts of an entire curve $f:\C \to A$ into a semi-abelian variety
$A$, provided that $f$ is of finite order.
On the other hand, by giving an example of $f$ into
a three dimensional abelian variety we show that this product
decomposition
does not hold in general; 
there was a gap in the proofs of \cite{N98}, Proposition~1.8 (ii) and of
\cite{SY96}, Theorem~2.2.
\end{abstract}

\section{Introduction}

Let $A$ be a semi-abelian variety
(i.e.~an extension of an abelian variety
by a product of copies of the multiplicative group), 
$f:\C\to A$ a holomorphic map
and $J_k(f):\C\to J_k(A)$ its $k$-th jet lift.
Since $A$ has a trivial tangent bundle, all the jet bundles
are trivial and $J_k(A)$ decomposes as a direct product
$A\times\C^{nk}$ with $n=\dim A$ in a canonical way.

In this note we use Zariski topology in terms of algebraic subsets.
A map from a complex space into a complex algebraic variety
is said to be Zariski (resp.\ non)degenerate if its image
is (resp.\ not) Zariski dense.

Now let $X_k$ denote the Zariski closure of
the image $J_k(f)(\C)$ in $J_k(A)$. Then it was claimed
in \cite{N98}, Proposition~ 1.8 (ii) and in
\cite{SY96}, Theorem~2.2 for the case of an abelian variety $A$
that $X_k$ is a direct product of a translate
of a semi-abelian subvariety of $A$ with a subvariety of $\C^{nk}$.
Unfortunately, there was a gap in these proofs, as pointed out by 
P.~Vojta in a private communication.

In this note we prove this assertion for the case where $f$ is of
finite order and also for the case where $A$ is a simple
abelian variety (see Theorem \ref{product}).
On the other hand, by giving an example we show 
that the assertion is not true in general (see Proposition \ref{counter}).
\section{The finite order case}

Basically, we follow the arguments of \cite{N98}. The key idea
is to consider jets of jets. For the better understanding
the argument is split in a sequence of lemmas.

\begin{lemma}
\label{hol}
Let $X$ be an irreducible complex algebraic variety and 
$p\in X\setminus\Sing(X)$  a non-singular point.
Then there exists a Zariski nondegenerate holomorphic map $f$ from 
the unit disk $\Delta$ to $X$ with $f(0)=p$.
\end{lemma}

\begin{proof}
Let $m=\dim X$. Let $\alpha_1,\ldots,\alpha_m$ be $\Q$-linearly
independent positive real numbers. 
Let $\Delta_2=\{z:|z|<2\}$ and fix an open embedding
$\iota :\Delta_2^m \hookrightarrow X$
such that $p\in \iota(\Delta_2^m)$.
Recall that there exists an isomorphicm
 $\zeta:\Delta \cong H^+=\{z:\Im(z)>0\}$.
We define a map $\phi:H^+\cup\R\to\Delta_2^m \subset X$ by
\[
z\mapsto \left( e^{\alpha_1iz},\ldots,e^{\alpha_miz}\right)
\]
Because the real numbers $\alpha_i$ are $\Q$-linearly independent,
$\phi(\R)$ is dense in the real torus $\left(S^1\right)^m$.
Hence $(S^1)^m$ is contained in the closure of $\phi(H^+)$.
Next we choose a holomorphic automorphism $\psi$ of $\Delta_2^m$
with $\psi(\phi(\zeta(0)))=p$ and define $f=\psi\circ\phi\circ\zeta$.

Because $(S^1)^m$ and therefore also $\psi\left((S^1)^m\right)$
is totally real and of real dimension $m$,
it follows that no closed analytic subset of $X$ except $X$
itself can contain $\iota(\psi(\phi(H^+)))$.
Thus $f(\Delta)$ must be Zariski dense.
\end{proof}

\begin{remark}\label{c-dense}\rm 
Regarded as a map from $\C$ to $\C^n$ the map $\phi$ defined in the
above proof is an example for a 
holomorphic map from $\C$ to $\C^n$ for which the image $\phi(\C)$ is
not contained in any proper analytic subset of $\C^n$.
\end{remark}

Let $A$ be a semi-abelian variety and $\C^m$ the complex affine space.
Let $X \subset A\times\C^m$ be an irreducible algebraic subvariety.
We define the stabilizer group by
\[
\mathrm{Stab}_A(X)=\{a\in A: (x+a,y)\in X\, \hbox{ for }
 \forall (x,y)\in X\subset A\times\C^m\}^0,
\]
where $\{\cdot\}^0$ stands for the identity component.
Then $\mathrm{Stab}_A(X)$ is a connected
closed algebraic subgroup of $A$.

For $l\in\N$, let $J_l(X)$ be the $l$-th jet space of $X$ and
let $\rho_l$ denote the natural projection
\[
\rho_l:J_l(A\times\C^m)\cong A\times \C^{nl+m(l+1)} 
\to J_l(A\times\C^m)/A \cong \C^{nl+m(l+1)} .
\]

\begin{lemma}
\label{stab}
Let the notation be as above.
Assume that $\mathrm{Stab}_A(X)=\{0\}$.
Then for every sufficiently large $l$ the differential
of $\rho_l|_{J_l(X)}$ restricted
to the jet space $J_l(X)$ has maximal rank at general points of
$J_l(X\setminus\Sing(X))$.
\end{lemma}

\begin{proof}
The proof is similar to that of Lemma (1.2) in \cite{N98}
(cf.\ \cite{NO}, Proof of Lemma (6.3.10), too).

For an arbitrarily fixed point $y_0 \in X\setminus\Sing(X)$ there is a
Zariski nondegenerate holomorphic map $f:\Delta \to X$ with
$f(0)=y_0$ by Lemma \ref{hol}.
Let $J_l(f):\Delta \to J_l(X)$ be the $l$-th jet lift, and set
$$
y_l=J_l(f)(0) \in J_l(X),\qquad l=0,1,2,\ldots .
$$
We are going to show that
there is a number $l_0$ satisfying the condition:
For every $l \geq l_0$, the differential at $y_l$
$$
d(\rho_l|_{J_l(X)})_{y_l}: \mathbf{T}_{y_l}(J_l(X)) \to
\mathbf{T}_{\rho_l(y_l)}(\C^{nl+m(l+1)}).
$$
is injective.

Let $x_0\in A$ be the image point of $y_l$ by the natural
projection, $J_l(X) \to A$.
First we note that
$$
\mathbf{T}_{y_l}(J_l(X)) \subset 
\mathbf{T}_{y_l}(A\times \C^{nl+m(l+1)})
\cong \mathbf{T}_{x_0}(A) \oplus
\C^{nl+m(l+1)}.
$$
Because of the definition of $\rho_l|_{J_l(X)}$,
we have the kernel,
$$
\ker d(\rho_l|_{J_l(X)})_{y_l} \subset \mathbf{T}_{x_0}(A).
$$
Since
$$
\ker d(\rho_l|_{J_l(X)})_{y_l} \supset
 \ker d(\rho_{l+1}|_{J_{l+1}(X)})_{y_{l+1}},
$$
there is a number $l_0$ such that
$$
\ker d(\rho_l|_{J_l(X)})_{y_l}=
\ker d(\rho_{l_0}|_{J_{l_0}(X)})_{y_{l_0}}
\hbox{ for all } l\geqq l_0.
$$
Suppose $\ker d(\rho_{l_0}|_{J_{l_0}(X)})_{y_{l_0}}\not=\{0\}$.
Take a vector
$v \in \ker d(\rho_{l_0}|_{J_{l_0}(X)})_{y_{l_0}}\setminus\{0\}$.
Let $h \in \mathcal{I}_{y_0}(X)$ be a holomorphic function germ
in the ideal sheaf of $X$ at $y_0$.
Then $v$, considered as a vector field on $A$, satisfies
$$
\left. \frac{d^l}{dz^l} \right|_{z=0} vh(f(z))=0,\qquad l=0,1,2,\ldots .
$$
Therefore, $vh(f(z))\equiv 0$ in a neighborhood of $0$
(cf.\ \cite{NO}, (6.3.12)). Since $h$ was
chosen arbitrarily,
we obtain $vh(f(z))\equiv 0$ near $0$ for all
$h\in \mathcal{I}_{y_0}(X)$. Hence $v$
is tangent to $X$ at every point of $f(\Delta)$.
Since $f (\Delta)$ is Zariski dense in $X$,
$v$ is everywhere tangent to $X$, so that
$X$ is invariant by the action of the one parameter subgroup
generated by $v$.
This contradicts the assumption.
\end{proof}

By Lemma \ref{stab} we immediately have the following.

\begin{lemma}
\label{alg}
Let $A$ be a semi-abelian variety,
and $g:\C\to A\times \C^m$ an entire curve.
Let $X$ be the Zariski closure of $g(\C)$. Assume that
$\mathrm{Stab}_A (X)=\{0\}$.
Let $l$ be as in Lemma \ref{stab}.
Then for every rational function $\phi$ on $X$
the induced meromorphic function $\phi \circ g$
on $\C$ is algebraic over the field extension of
$\C$ generated by all the components of
$\rho_l\circ J_l(g)$.
\end{lemma}

\begin{theorem}
\label{product}
Let $A$ be a semi-abelian variety of dimension $n$,
and $f:\C\to A$ an entire curve.
Let $X_k$ be the Zariski closure of $J_k(f)(\C)$
in $J_k(A) \cong A\times \C^{nk}$.
\begin{enumerate}
\item
Assume that $f$ is of finite order.
Then there exist a semi-abelian subvariety $B\subset A$,
$a \in A$, and
a subvariety $W_k\subset \C^{nk}$
such that $X_k=(B+a)\times W_k$.
\item
If $A$ is a simple abelian variety, then there is
a subvariety $W_k\subset \C^{nk}$
such that $X_k=A\times W_k$.
\end{enumerate}
\end{theorem}

\begin{proof}
Let $B=\st_A(X_k)$
and apply the preceding
Lemma \ref{alg} to the holomorphic map
$F:\C\to (A/B)\times\C^{nk} \cong J_k(A)/B$
induced by $J_k(f):\C \to J_k(A)\cong A \times \C^{nk}$.
Let $F_1:\C\to A/B$ denote its first component.
Then $F_1$ is the composition of $f:\C\to A$ with the natural
projection from $A$ to $A/B$.

Assume that $f$ is of finite order.
Let $\tilde f:\C\to\C^n$ denote the lift of $f$ to the universal covering
of $A$.
Then the lift $\tilde f$ is a polynomial
map, and hence $\rho_l\circ J_l(F)$ are polynomial map, too
(see, e.g., \cite{NWY02}).
The preceding Lemma \ref{alg} now yields that
for every rational function $\phi$
on the semi-abelian variety $A/B$ the induced meromorphic function
$\phi \circ F_1$ is algebraic over the rational function
field $\C(z)$ in $z$.
It follows that Nevanlinna's order function fulfills 
$T_{\phi \circ F_1}(r)=O(\log r)$
(cf., e.g., \cite{NO}, Proposition (5.3.14) and Lemma (6.1.5)).
This implies that $\phi\circ F_1$ is a rational function.
Since this is true for every rational function $\phi$ on $A/B$,
it follows that $F_1:\C\to A/B$ is an algebraic morphism.
But every algebraic morphism from $\C$ to
a semi-abelian variety is constant. Thus the projection map from
$J_k(A)$ to $A/B$ maps $X_k$ to a point. This implies
that $X_k=(B+a)\times W_k$ for some subvariety $W_k\subset J_k(A)/A$
and $a \in A$.

In general, the above argument at least
implies $\dim\st(X_k)>0$ (cf.\ \cite{N98}, Lemma (1.2)).
If $A$ is a simple abelian variety, then
$\st(X_k)=A$, so that the required assertion follows.
\end{proof}

\begin{remarks}
\null\ 
\rm
\begin{enumerate}
\item
Using induction, one can deduce from
  $\dim \mathrm{Stab}_A(X_0)>0$
that for every such holomorphic
map $f$ the Zariski closure of the image $f(\C)$ in $A$ is a
translate of a semi-abelian subvariety.
This fact (``logarithmic Bloch-Ochiai theorem'') was first
proved in \cite{N81}.
\item
However, for $k>0$ and $0<\dim B <\dim A$
the quotient $J_k(A)/B$ is
strictly larger than $J_k(A/B)$.
Therefore one cannot use induction to prove that $X_k$ has a 
product decomposition as in the case where $f$ is of finite order.
Basically, this is the gap in the incomplete proofs mentioned in the
introduction.
This gap cannot be filled:
As we will see in the next section there is a counter-example
in the case where $f$ is not of finite order.
\end{enumerate}
\end{remarks}

\section{The counter-example}

\begin{proposition}
\label{counter}
There exist a three-dimensional complex abelian variety $A$ and a 
holomorphic map $f:\C\to A$ with the following properties:
\begin{enumerate}
\item
The image $f(\C)$ is Zariski dense in $A$;
\item
Let $X$ denote the Zariski closure of the image of the jet lift
$J_1(f):\C\to J_1(A)\cong A\times\C^3$.
Then $X$ is {\em not} a direct product inside $J_1(A)$, i.e.,
there does not exist a subvariety $W\subset \C^3$
such that $X=A \times W$.
\end{enumerate}
\end{proposition}
\begin{proof}
Let $C$ be a nonsingular cubic curve in $\P_2(\C)$. Then $C$ is an
elliptic curve, i.e., a one-dimensional compact complex torus.
For a point $x=[x_0:x_1:x_2]\in\P_2(\C)$ we define
$\ortho{x}=\{[z_0:z_1:z_2]:\sum_i z_ix_i=0\}$.
We define a surface $\bar S\subset C\times\P_2(\C)$ by
\[
\bar S=\{(x,z)\in \P_2(\C)\times\P_2(\C): z\in \ortho{x}\text{ and }x\in C
\}
\]
Note that the projection $C\times\P_2(\C)\to C$ induces
the structure of a holomorphic fiber bundle on $\bar S$, with
fiber $\P_1(\C)$ and base $C$. 

Fix $p\in\P_2(\C)\setminus C$.
Now for every $x\in C$, the sets $\ortho{p}$ and $\ortho{x}$
are two different lines in $\P_2(\C)$. In particular
$\ortho{x}\cap\ortho{p}$ contains exactly one point.
Define $S=\bar S\setminus (C\times \ortho{p})$.
Then $S$ is a closed algebraic subvariety of $C\times\C^2\cong C
\times\left(\P_2(\C)\setminus\ortho{p}\right)$.
Let $i:S\to C\times\C^2$ denote this embedding.
On the other hand, by the projection on the first factor,
$S$ can be realized as a holomorphic fiber bundle with fiber $\C$ and
base $C$. 
Let $\pi:\tilde S\to S$ be the universal covering.
Then $\tilde S\stackrel{\varphi}\cong\C^2$, because $\C$ is the universal
covering of $C$ and every $\C$-bundle over $\C$ is trivial.

Let $g: \C \to \C^2$ be a holomorphic map such that no proper
analytic subset of $\C^2$ contains $g(\C)$ (cf.\ Remark \ref{c-dense}).
Let $G=(G_1, G_2):\C\to C\times\C^2$ be the holomorphic map given by
the compositions
\[
\C\stackrel{g}{\longrightarrow}\C^2 \stackrel{\varphi}\cong\tilde S
\stackrel{\pi}{\longrightarrow}S\stackrel{i}{\longrightarrow}C\times\C^2 .
\]
Note that $i(S)$ is the Zariski closure of $G(\C)$ in $C\times\C^2$.

Choose a lattice $\Lambda$ in $\C^2$ such that $B=\C^2/\Lambda$
is a simple complex abelian variety.
Denote the natural projection by $\tau:\C^2\to B$.
Let $A=C\times B$.
Now we are in a position to define the desired entire curve
 $f:\C\to A$:
Let $f=(f_1,f_2):\C\to C\times B$
with $f_1=G_1$ and
\[
f_2(x)=\tau\left( \int_0^x G_2(t) dt \right ).
\]

We have to verify that $f$ has the desired properties.
By construction $f(\C)$ is not contained in (a translate of)
one of the two factors
of the product $C\times B$.
On the other hand,
the Zariski closure of $f(\C)$ must be (a translate of) a subtorus and
$B$ is simple. Therefore $f(\C)$ is Zariski dense in $A$.

Next we consider the first jet lift
\[
J_1(f):\C\to J_1(A)\cong A\times\C^3 \cong
 C \times B \times \C \times \C^2 .
\]
By construction we have
\[
J_1(f)=(G_1,f_2,G_1',G_2) .
\]
Let $q:(x_1,x_2;v_1,v_2)\mapsto (x_1,v_2)$ denote the projection on the
first and fourth factor.
By our construction the closure of $q(J_1(f)(\C))$
in $C\times\C^2$ coincides with $i(S)$. In particular,
it is not a direct product of subvarieties of the two factors.
It follows that the Zariski
closure of $J_1(f)$ in $J_1(A)=A\times\C^3$ 
cannot be a direct product of a subvariety of $A$ with a subvariety
of $\C^3$.
\end{proof}

\begin{remark}\rm
The projection map $p_2$ from $i(S)\subset C\times \C^2$
onto the second factor $\C^2$ maps $i(S)$ surjectively onto $\C^2$.
As a result, for a holomorphic map $G:\C\to i(S)$ with Zariski dense 
image in $i(S)$ the induced map $G_2=p_2\circ G:\C\to\C^2$ has
a Zariski dense image.
Therefore in the above construction $G_2:\C\to\C^2$ 
cannot be a polynomial map. 
As a consequence $f_2$ and therefore $f:\C\to A$
is not a holomorphic map of finite order.

Note also that $\st_A(X)=\{0\}\times B$.
\end{remark}

\section{Consequences}
It should be noted that in the proof of the Main Theorem
of \cite{N98}, Proposition (1.8) (ii) was not used. Hence it is
valid without any change.

Similarly, the product decomposition of $X_k$ is not
essential in any of the proofs of \cite{NWY02} as in
\cite{NWY00}, where the product decomposition was not used.

\end{document}